\documentclass[11pt]{article}
\usepackage{graphicx}
\usepackage{amsmath}
\usepackage{epsfig}  % to include an Encapsulated PostScript graphic
\usepackage{verbatim}  
\usepackage{tabularx}
\usepackage{epic, eepic}
\usepackage{xcolor}
\usepackage{lineno}
\usepackage{amssymb,latexsym,amscd,amsmath,amsfonts,enumerate,supertabular}
\usepackage{subfigure}
\usepackage{longtable}
\newtheorem{Theorem}{Theorem}[section] 

\newtheorem{lemma}[Theorem]{Lemma} 
\newtheorem{Proposition}[Theorem]{Proposition}

\begin{document}

\title{{\bf Existence results and iterative method for solving a  fourth order nonlinear integro-differential equation }}

%\titlerunning{ Simple numerical method for solving a fully third order nonlinear BVP}        % if too long for running head
\author{ Dang Quang A$^{\text a}$,  Dang Quang Long$^{\text b}$\\
$^{\text a}$ {\it\small Center for Informatics and Computing, VAST}\\
{\it\small 18 Hoang Quoc Viet, Cau Giay, Hanoi, Vietnam}\\
{\small Email: dangquanga@cic.vast.vn}\\
$^{\text b}$ {\it\small Institute of Information Technology, VAST,}\\
{\it\small 18 Hoang Quoc Viet, Cau Giay, Hanoi, Vietnam}\\
{\small Email: dqlong88@gmail.com}}
\date{ }
\date{}          % Enter your date or \today between curly braces
\maketitle
\begin{abstract}
In this paper we consider a class of  fourth order nonlinear integro-differential equations with Navier boundary conditions.  By the reduction of the problem to operator equation we establish the existence and uniqueness of solution and construct a numerical method for solving it. We prove that the method is of second order accuracy and obtain an estimate for total error.  Some examples demonstrate the validity of the obtained theoretical results and the efficiency of the numerical method.
\end{abstract}
{\small
}
\section{Introduction}
Integro-differential equations are the mathematical models of many phenomena in physics, biology, hydromechanics, chemistry, etc. In general, it is impossible to find the exact solutions of the problems involving these equations, especially when they are nonlinear. Therefore, many analytical approximate  methods and numerical methods are developed for  these equations (see, e.g., \cite{Aruc,Chen,Dasc,Lake,Sing,Swei,Tahe,Yula,Wang,Zhua}).

Below, we mention some works concerning the solution methods for integro-differential equations. First, it is worthy to mention the recent work of Tahernezhad and Jalilian  in 2020 \cite{Tahe}. In this work the authors consider the second order linear problem

\begin{align*}
u''(x)&+p(x)u'(x)+q(x)u(x)=f(x)+ \int_a^b k(x,t)u(t) dt , \quad a<x<b,\\
u(a)&=\alpha , \; u(b)=\beta, \end{align*}
where $p(x), q(x), k(x,t)$ are sufficiently smooth functions. \par 
Using the non-polynomial spline functions, namely, the exponential spline functions, the authors constructed the numerical solution of the problem and proved that the error of the approximate solution is $O(h^2)$, where $h$ is the grid size on $[a, b]$. Before \cite{Tahe} there is an interesting works of Chen et al. \cite{Chen,Chen1}, where the authors used  multiscale Galerkin method for constructing the approximate solutions of the above second order problem, for which the computed convergence rate is two. \par
Besides the researches evolving the second order integro-differential equations, recently many authors are interested in fourth order integro-differential equations due to their wide applications. We first mention the work of Singh and Wazwaz \cite{Sing}. In this work the authors developed a technique based on the Adomian decomposition method with the Green function for constructing  a series solution of the nonlinear Voltera equation associated with the Dirichlet boundary conditions
\begin{align}
y^{(4)}(x)&=g(x)+ \int_0^x k(x,t)f(y(t)) dt , \quad 0<x<b, \label{eq:intr1} \\
y(0)&=\alpha_1 , \; y'(0)=\alpha_2, \; y(b)=\alpha_3 , \; y'(b)=\alpha_4. \label{eq:intr2}
\end{align}
Under some conditions it was proved that the series solution converges as a geometric progression.\par
For the linear Fredholm IDE \cite{Aruc}
\begin{align*}
y^{(4)}(x)+\alpha y''(x) +\beta y (x)- \int_a^b K(x,t)y(t) dt = f(x)  , \quad a<x<b,
%y(a)&=\alpha_1 , \; y'(a)=\alpha_2, \; y(b)=\alpha_3 , \; y'(b)=\alpha_4. 
\end{align*}
with the Dirichlet boundary conditions \eqref{eq:intr2}, the difference method and the trapezium rule are used to design the corresponding linear system of algebraic equations. A new variant called Modified Arithmetic Mean iterative method is proposed for solving the latter system, but the error estimate of the method is not obtained.\par
%In \cite{Swei} the 
The boundary value problem for the nonlinear IDE
\begin{align*}
y^{(4)}(x)&-\varepsilon y''(x)-\frac{2}{\pi}\Big ( \int_0^\pi |y'(t)|^2 dt  \Big )y''(x) =p(x) , \quad 0<x<\pi,  \\
y(0)&=0 , \; y(\pi)=0, \; y''(0)=0 , \; y''(\pi)=0 
\end{align*}
was considered in \cite{Dang-Vu,Zhua}, where the authors constructed approximate solutions by the iterative and spectral methods, respectively. Recently, Dang and Nguyen \cite{Dang4} studied the existence and uniqueness of solution and constructed iterative method for finding the solution for the IDE
\begin{align*}
u^{(4)}(x)&-M\Big ( \int_0^L |u'(t)|^2 dt  \Big )u''(x) = f(x,u,u',u'',u''') , \quad 0<x<L,  \\
u(0)&=0 , \; u(L)=0, \; u''(0)=0 , \; u''(L)=0, 
\end{align*}
where $M$ is a continuous non-negative function. \par 
Very recently, Wang \cite{Wang} considered the problem
\begin{equation*}
\begin{split}
y^{(4)}(x)& =f(x,y(x), \int_0^1 k(x,t) y(t) dt ),\\
y(0)&=0 , \; y(1)=0, \; y''(0)=0 , \; y''(1)=0.
\end{split}
\end{equation*}
Using the monotone method, the author constructed the sequences of functions, which converge to the extremal solutions.\par
Motivated by the above facts, in this paper we consider an extension of the above problem, namely, the problem 
\begin{equation}\label{eq:O3}
\begin{split}
u^{(4)}(x)& =f(x,u(x),u'(x), \int_0^1 k(x,t) u(t) dt ),\\
u(0)&=0 , \; u(1)=0, \; u''(0)=0 , \; u''(1)=0,
\end{split}
\end{equation}
where the function $f(x,u,v,z)$ and $k(x,t)$ are assumed to be continuous.
Using the method developed in our previous papers \cite{Dang1,Dang2,Dang4,Dang5,Dang6,Dang7} we establish the existence and uniqueness of solution and propose an iterative method at both continuous and discrete levels for finding the solution. The second order convergence of the method is proved. The theoretical results are illustrated on some examples.

\section{Existence results}
Using the methodology in \cite{Dang1,Dang2,Dang4,Dang5,Dang6,Dang7} we introduce the operator $A$ defined in the 
space of continuous functions $C[0, 1]$   by the formula
\begin{equation}\label{eq:A1}
(A\varphi ) (x)= f(x,u(x),u'(x),\int_0^1 k(x,t) u(t) dt ),
\end{equation}
where $u(x)$ is the solution of the boundary value problem
\begin{equation}\label{eq:A2}
\begin{split}
u^{\prime\prime\prime\prime}&=\varphi (x), \; 0<x<1,\\
u(0)&=u^{\prime\prime}(0)=u(1)=u^{\prime\prime}(1)=0 .
\end{split}
\end{equation}
It is easy to verify the following lemma.
\begin{lemma}
If the function $\varphi$ is a fixed point of the operator $A$, i.e., $\varphi$ is the solution of the operator equation 
\begin{equation}\label{eq:A3}
A\varphi = \varphi ,
\end{equation} 
where $A$ is defined by \eqref{eq:A1}-\eqref{eq:A2} 
then the function $u(x)$ determined from the BVP \eqref{eq:A2} is a solution of the BVP \eqref{eq:O3}. Conversely, if the function $u(x)$ is the solution of the BVP \eqref{eq:O3} then the function 
\begin{equation*}
\varphi (x)= f(x,u(x),u'(x),\int_0^1 k(x,t) u(t) dt )
\end{equation*}
satisfies the operator equation \eqref{eq:A3}.
\end{lemma}
Due to the above lemma we shall study the original BVP \eqref{eq:O3} via the operator equation \eqref{eq:A3}. Before doing this we notice that  the BVP \eqref{eq:A2} has a unique solution representable in the form
\begin{equation}\label{eq:A4}
u(x)=\int_0^1 G_0(x,s)\varphi(s)ds, \quad 0<x<1,
\end{equation}
where 
\begin{align}\label{eq:A5}
G_0(x,s)=\frac{1}{6}
\begin{cases}
s(x-1)(x^2-x+s^2), \quad & 0\leq s\leq x\leq 1 \\
x(s-1)(s^2-s+x^2), & 0\leq x\leq s\leq 1
\end{cases}
\end{align}
is the Green function of the operator $u''''(x)=0$ associated with the homogeneous boundary conditions 
$u(0)=u''(0)= u(1)=u''(1)=0$.\\
Differentiating both sides of \eqref{eq:A4} gives
\begin{align}
u'(x) & =\int_0^1 G_1(x,s)\varphi(s)ds, \label{eq:A6}
\end{align}
where
\begin{equation}\label{eq:A7}
G_1(x,s)  =\frac{1}{6}
\begin{cases}
s(3x^2-6x+s^2+2), \quad & 0\leq s\leq x\leq 1, \\
(s-1)(3x^2-2s+s^2),& 0\leq x\leq s\leq 1 .
\end{cases}   
\end{equation}
Set 
\begin{align}\label{eq:A8}
\begin{split}
M_0 & =\max_{0\leq x\leq 1} \int_0^1 |G_0(x,s)|ds ,\\
M_1 & =\max_{0\leq x\leq 1} \int_0^1 |G_1(x,s)|ds , \\
M_2 & =\max_{0\leq x\leq 1} \int_0^1 |k(x,s)|ds
\end{split}
\end{align}  
It is easy to obtain 
\begin{equation}\label{eq:M}
 M_0=   \frac{5}{384} , M_1= \frac{1}{24}.
\end{equation}

Now for any positive number $M$, we define the domain
\begin{align}\label{eq:A9}
\begin{split}
\mathcal{D}_M = \{(x,u,v,z)\; | \; &0\leq x\leq 1,\; |u|\leq M_0M,\\
 &|v|\leq M_1M,\;  |z|\leq M_0M_2M \}.
\end{split}
\end{align}

\begin{Theorem}[Existence and uniqueness]\label{thm1}
Suppose that the function $k(x,t)$ is continuous in the square $[0,1] \times [0,1]$ and there exist numbers $M>0$, $L_0,L_1,L_2\geq 0$ such that:
\begin{description}
\item [(i)]The function $f(x,u,v,z)$ is continuous in the domain $\mathcal{D}_M$ and $|f(x,u,v,z)|\leq M$, $ \forall (x,u,v,z)\in \mathcal{D}_M$.
\item [(ii)]$|f(x_2,u_2,v_2,z_2)-f(x_1,u_1,v_1,z_1)| \leq L_0 |u_2-u_1|+ L_1 |v_2-v_1|+ L_2 |z_2-z_1|, \; \forall (x_i,u_i,v_i,z_i)\in \mathcal{D}_M, \; i=1,2.$
\item [(iii)]$q = L_0M_0+L_1M_1+L_2M_0M_2 <1$
\end{description}
%\end{enumerate}
Then the problem \eqref{eq:O3} has a unique solution $u\in C^4[0,1]$ satisfying $|u(x)| \leq M_0M,$  $ |u'(x)| \leq M_1M$ for any $0 \le x \le 1.$
\end{Theorem}
%\begin{Proof}
\noindent {\bf Proof.}
Under the assumptions of the theorem we shall prove that the operator $A$ is a contraction mapping in the closed ball $B[O,M]$. Then the operator equation \eqref{eq:A3} has a unique solution $u \in C^{(4)}[0,1]$ and this implies the existence and uniqueness of solution of the BVP \eqref{eq:O3}.\par 
Indeed, take $\varphi \in B[O,M]$. Then the problem \eqref{eq:A2} has a unique solution of the form \eqref{eq:A4}. From there and \eqref{eq:A8} we obtain $|u(x)| \le M_0\| \varphi \|$ for all $x\in [0.1]$. Analogously, we have $\| u'(x)\| \le M_1 \|\varphi \|$ for all $x\in [0.1]$. Denote by $K$ the integral operator defined by
$$(Ku)(x) = \int_0^1 k(x,t) u(t) dt.$$
Then from the last equation in \eqref{eq:A8} we have the estimate $| (Ku)(x)| \le M_2 \|\varphi\|, \; x\in [0,1]$. Thus, if $\varphi \in B[O,M]$, i.e., $\|\varphi\| \le M$ then for any $x \in [0,1]$ we have
$$ |u(x)| \le M_0M, \; |u'(x)| \le M_1M, \; |(Ku) (x)|\le M_0M_2M.
$$
Therefore, $(x, u(x),u'(x), (Ku) (x))\in \mathcal{D}_M$. By the assumption $(i)$ there is
$$ |f(x, u(x),u'(x), (Ku) (x))| \le M \quad \forall x\in [0,1].
$$
Hence, $|(A\varphi)(x)| \le M \forall x\in [0,1]$ and $\|A\varphi \| \le M$. It means that $A$ maps $B[O,M]$ into itself.\par 
Next, take $\varphi _1, \varphi_2 \in B[O,M]$. Using the assumption $(ii)$ and $(iii)$ it is easy to obtain
\begin{align*}
\|A\varphi _2 - A\varphi _1\| \le (L_0M_0+L_1M_1+L_2M_0M_2)\|\varphi _2 -\varphi_1 \| =q \|\varphi _2 -\varphi_1 \|.
\end{align*}
Since $q<1$ the operator $A$ is a contraction in $B[O,M]$. This completes the proof of the theorem. $\blacksquare$

%\end{Proof}
Now, in order to study positive solutions of the BVP \eqref{eq:O3} we introduce the domain
\begin{align}\label{eq:A10}
\begin{split}
\mathcal{D}_M^+ = \{(x,u,v,z)\; | \; &0\leq x\leq 1,\; 0\le u\leq M_0M,\\
 &|v|\leq M_1M,\;  |z|\leq M_0M_2M \}.
\end{split}
\end{align}
and denote
$$ S_M = \{\varphi \in C[0,1], \|\varphi \| \le M   \}
$$
\begin{Theorem}[Positivity of solution]\label{thm2}
Suppose that the function $k(x,t)$ is continuous in the square $[0,1] \times [0,1]$ and there exist numbers $M>0, L_0,L_1,L_2\geq 0$ such that:
\begin{description}
\item [(i)]The function $f(x,u,v,z)$ is continuous in the domain $\mathcal{D}_M^+$ and $0\le f(x,u,v,z)\leq M, \; \forall (x,u,v,z)\in \mathcal{D}_M^+$ and $f(x,0,0,0) \not \equiv 0$ 
\item [(ii)]$|f(x_2,u_2,v_2,z_2)-f(x_1,u_1,v_1,z_1)| \leq L_0 |u_2-u_1|+ L_1 |v_2-v_1|+ L_2 |z_2-z_1|, \; \forall (x_i,u_i,v_i,z_i)\in \mathcal{D}_M^+, \; i=1,2.$
\item [(iii)]$q = L_0M_0+L_1M_1+L_2M_2 <1$
\end{description}
Then the problem \eqref{eq:O3} has a unique positive solution $u\in C^4[0,1]$ satisfying $0 \le u(x) \leq M_0M,$  $ |u'(x)| \leq M_1M$ for any $0 \le x \le 1.$
\end{Theorem}
%\begin{Proof}
\noindent {\bf Proof.}
Similarly to the proof of Theorem \ref{thm1}, where instead of $\mathcal{D}_M$ and $B[O; M]$
there stand $\mathcal{D}_M^+$ and $S_M$, we conclude that the problem  has a nonnegative solution. Due to the condition $f(x,0,0,0) \not \equiv 0$, this solution must be positive.  $\blacksquare$

%\end{Proof}

\section{Iterative method }\label{IterMeth}
Consider the following iterative method
\begin{enumerate}
\item Given 
\begin{equation}\label{iter1c}
\varphi_0(x)=f(x,0,0,0),
\end{equation}
\item Knowing $\varphi_m(x)$  $(m=0,1,...)$ compute
\begin{equation}\label{iter2c}
\begin{split}
u_m(x) &= \int_0^1   G_0(x,t)\varphi_m(t)dt  ,\\
v_m(x) &= \int_0^1 G_1(x,t)\varphi_m(t)dt ,\\
z_m(x) &= \int_0^1 k(x,t)u_m(t) dt 
\end{split}
\end{equation}
\item Update
\begin{equation}\label{iter3c}
\begin{split}
\varphi_{m+1}(x) &= f(x,u_m(x),v_m(x),z_m(x)).
\end{split}
\end{equation}
\end{enumerate}
This iterative method indeed is the successive iterative method for finding the fixed point of operator $A$. Therefore, it converges with the rate of geometric progression and there holds the estimate
$$\|\varphi_m -\varphi\| \leq \frac{q^m}{1-q}\|\varphi_1 -\varphi_0\| = p_m d,$$
where  $\varphi$ is the fixed point of the operator $A$ and
\begin{equation}\label{eq:pd}
p_m=\tfrac{q^m}{1-q},\; d=\|\varphi_1 -\varphi_0\|.
\end{equation}

These estimates imply the following result of the convergence of the iterative method \eqref{iter1c}-\eqref{iter3c}.
\begin{Theorem}\label{thm3}
The iterative method \eqref{iter1c}-\eqref{iter3c} converges and for the approximate solution $u_k(t)$ there hold estimates
\begin{align*}
\|u_m-u\| &\leq M_0p_md, \; \|u'_m-u'\| \leq M_1p_md, 
\end{align*}
where $u$ is the exact solution of the problem \eqref{eq:O3}, $p_m$ and $d$ are defined by \eqref{eq:pd}.
\end{Theorem}

%Now, in order to realize the iterative method at continuous level consider an iterative method at discrete level. 

%\section{Discrete iterative method}
To numerically realize the above iterative method we construct  corresponding discrete iterative method. For this purpose cover the interval $[0, 1]$   by the uniform grid $\bar{\omega}_h=\{x_i=ih, \; h=1/N, i=0,1,...,N  \}$ and denote by $\Phi_m(x), U_m(x),  V_m(x), Z_m(x)$ the grid functions, which are defined on the grid $\bar{\omega}_h$ and approximate the functions $\varphi_m (x), u_m(x),  v_m(x), z_m(x)$ on this grid. \par
Consider now the following discrete iterative method.
\begin{enumerate}
\item Given 
\begin{equation}\label{iter1d}
\Phi_0(x_i)=f(x_i,0,0,0),\ i=0,...,N; 
\end{equation}
\item Knowing $\Phi_m(x_i), \; i=0,...,N $  compute approximately the definite integrals \eqref{iter2c} by trapezium formulas
\begin{equation}\label{iter2d}
\begin{split}
U_m(x_i) &= \sum _{j=0}^N h\rho_j G_0(x_i,x_j)\Phi_m(x_j) ,\\
V_m(x_i) &= \sum _{j=0}^N h\rho_j G_1(x_i,x_j)\Phi_m(x_j) ,\\
Z_m(x_i) &= \sum _{j=0}^N h\rho_j k(x_i,x_j)U_m(x_j) ,\;  i=0,...,N,
\end{split}
\end{equation}
\noindent where $\rho_j$ is the weight of the trapezium formula
\begin{equation*}
\rho_j = 
\begin{cases}
1/2,\; j=0,N\\
1, \; j=1,2,...,N-1
\end{cases}
\end{equation*}

\item Update
\begin{equation}\label{iter3d}
\Phi_{m+1}(x_i) = f(x_i,U_m(x_i),V_m(x_i),Z_m(x_i)).
\end{equation}
\end{enumerate}
In order to get the error estimates for the approximate solution for $u(t)$ and its derivatives on the grid we need some following auxiliary results.  \\

\begin{Proposition}\label{prop1}
Assume that the function $f(t,u,v,z)$ has all continuous partial derivatives up to second order in the domain $\mathcal{D}_M$ and the kernel function $k(x,t)$ has all continuous partial derivatives up to second order in the  square $[0,1] \times [0,1]$. Then for the functions $\varphi_m(x), u_m(x),  v_m(x), z_m(x), m=0,1,...$ constructed by the iterative method \eqref{iter1c}-\eqref{iter3c} we have  
$\varphi_m(x) \in  C^2 [0, 1]$, $u_m(x) \in C^6 [0, 1]$, $ v_m(x) \in C^5 [0, 1], \; z_m(x) \in C^2 [0, 1]$
\end{Proposition}
%\begin{Proof}
\noindent {\bf Proof.}
We prove the proposition by induction. For $k=0,$ by the assumption on the function $f$ we have $\varphi_0(t) \in C^2[0, 1]$ since $\varphi_0(x)=f(x,0,0,0)$. Taking into account that
\begin{align*}
u_0(x) = \int_0^1   G_0(x,t)\varphi_0(t)dt
\end{align*}
the function $u_0(x)$ is the solution of the BVP
\begin{align*}
u_0^{(4)}(x) &= \varphi_0 (x),\quad  x\in (0,1),\\
u_0(0)&=u_0(1)= u_0''(0)=u_0''(1)=0.
\end{align*}
Therefore, $u_0(x) \in C^6 [0, 1]$.  It implies that $v_0(x) \in C^5 [0, 1]$ because $v_0(x)= u_0'(x)$. Since by assumptions $k(x,t)$ has all continuous derivatives up to second order, the function $z_0(x) =\int_0^1 k(x,t) u_0(t)dt$ belongs to $C^2[0,1]$.\par

Now suppose $\varphi_m(x) \in  C^2 [0, 1]$, $u_m(x) \in C^6 [0, 1]$, $ v_m(x) \in C^5 [0, 1], \; z_m(x) \in C^2 [0, 1]$.  Then, because
$\varphi_{m+1}(x) = f(x,u_m(x),v_m(x),z_m(x))$
and the functions $f$ by the assumption has continuous derivative in all variables up to order 2, it follows that $\varphi_{m+1}(x) \in C^2[0, 1]$. Repeating the same argument as for $\varphi_0(x)$ above we obtain that $u_{m+1}(x) \in C^6 [0, 1]$, $ v_{m+1}(x) \in C^5 [0, 1], \; z_{m+1}(x) \in C^2 [0, 1]$
Thus, the proposition is proved.  $\blacksquare$ \\
%\end{Proof}

\begin{Proposition}\label{prop2}

For any function $\varphi (x) \in C^2[0, 1]$ there hold the estimate
\begin{equation}
\int_0^1 G_n (x_i,t) \varphi (t) dt = \sum _{j=0}^N h\rho_j G_n(x_i,t_j)\varphi(t_j) +O(h^2), 
\quad (n=0,1), \label{eq:prop2}\\
\end{equation}
\end{Proposition}
%\begin{Proof}
\noindent {\bf Proof.} The above estimate is obvious in view of the error estimate of the compound trapezium formula because the function $G_n(x_i,t)\; (n=0,1)$  are continuous at $t_j$ and are polynomials in the intervals $[0,t_j]$ and $[t_j, 1]$.        $\blacksquare$
 
%\end{Proof} 

\begin{Proposition}\label{prop3}
Under the assumption of Proposition \ref{prop1}  for any $m=0,1,...$ there hold the estimates
\begin{equation}\label{eq:prop3a}
\|\Phi_m -\varphi_m  \|= O(h^2),\quad \|U_m -u_m  \|=O(h^2),
\end{equation}
\begin{equation}\label{eq:prop3b}
\begin{split}
  \|V_m-v_m  \|=O(h^2), \; \|Z_m -z_m  \|=O(h^2).
\end{split}
\end{equation}
 where $ \|.\|=  \|.\|_{C(\bar{\omega}_h)}$ is the max-norm of function on the grid $\bar{\omega}_h$.
\end{Proposition}
%\begin{Proof}
\noindent {\bf Proof.} We prove the proposition by induction. For $m=0$ we have immediately $\|\Phi_0 -\varphi_0  \|= 0$. Next, by the first equation in \eqref{iter2c} and Proposition \ref{prop2} we have
\begin{equation}
u_0(x_i)=\int_0^1 G_0 (x_i,t) \varphi_0 (t) dt = \sum _{j=0}^N h\rho_j G_0(x_i,t_j)\varphi_0(t_j)+O(h^2)
\end{equation}
for any $i=0,...,N$ .
On the other hand, in view of  the first equation in \eqref{iter2d} and \eqref{iter1d} we have
\begin{equation}
U_0(x_i)= \sum _{j=0}^N h\rho_j G_0(x_i,t_j)\Phi_0(t_j)
\end{equation}
Therefore, $|U_0(t_i)- u_0(t_i)|= O(h^2)$. Consequently, $\|U_0 -u_0  \|=O(h^2) $.\\
Similarly, we have 
\begin{equation}
\|V_0 -v_0  \|=O(h^2).
\end{equation}
By the trapezium formula we have
\begin{align*}
z_0(x_i)= \int_0^1 k(x_i,t)u_0(t) dt =  \sum _{j=0}^N h\rho_j k(x_i,t_j)u_0(t_j) +O(h^2),
\end{align*}
while by the third equation in \eqref{iter2d} we have
\begin{align*}
Z_0(x_i) = \sum _{j=0}^N h\rho_j k(x_i,t_j)U_0(t_j) ,\;  i=0,...,N.
\end{align*}
Therefore, 
\begin{align*} 
\Big | Z_0(x_i) - z_0(x_i)   \Big | &= \Big | \sum _{j=0} ^{N} h\rho_j k(x_i,t_j) (U_0(t_j) -u_0(t_j) )  \Big | +O(h^2) \\
&\le  \sum _{j=0} ^{N} h\rho_j |k(x_i,t_j)| |U_0(t_j) -u_0(t_j) |   +O(h^2) \\
&\le C h^2  \sum _{j=0} ^{N} h\rho_j |k(x_i,t_j)| + O(h^2) \\
&\le CC_1 \sum _{j=0} ^{N} h\rho_j + O(h^2) =  O(h^2)
\end{align*} 
because $|U_0(t_j) -u_0(t_j) | \le Ch^2, \; |k(x_i,t_j)| \le C_1$, where $C, C_1$ are some constants.

Now suppose that \eqref{eq:prop3a} and \eqref{eq:prop3b} are valid for $m \ge 0$. We shall show that these estimates are valid for $m+1$.
By the Lipschitz condition of the function $f$ and the estimates \eqref{eq:prop3a} and \eqref{eq:prop3b} it is easy to obtain the estimate $\|\Phi_{m+1} -\varphi_{m+1}  \|= O(h^2)$. 
Now from the first equation in \eqref{iter2c} by Proposition \ref{prop3} we have
\begin{equation*}
u_{m+1}(x_i)=\int_0^1 G_0 (x_i,t) \varphi_{m+1} (t) dt = \sum _{j=0}^N h\rho_j G_0(x_i,x_j)\varphi_{m+1}(x_j)+O(h^2)
\end{equation*}
On the other hand by the first formula in \eqref{iter2d} we have
\begin{equation*}
U_{m+1}(x_i) = \sum _{j=0}^N h\rho_j G_0(x_i,x_j)\Phi_{m+1}(x_j).
\end{equation*}
From this equality and the above estimates  we obtain the estimate
$$\|U_{m+1} -u_{m+1}  \|=O(h^2).
$$
Similarly, we obtain
$$\|V_{m+1} -v_{m+1}\|=O(h^2),\; \|Z_{m+1} -z_{k+1}\|=O(h^2).
$$
Thus, by induction we have proved the proposition.  $\blacksquare$ \par
%\end{Proof}
Now combining Proposition \ref{prop3} and Theorem \ref{thm1} results in the following theorem.
\begin{Theorem}\label{thm4}
For the approximate solution of the problem \eqref{eq:O3} obtained by the discrete iterative method on the uniform grid with grid size $h$ there hold the estimates
\begin{equation}\label{eqthm4}
%\begin{split}
\|U_k-u\| \leq M_0 p_kd +O(h^2), \; \|V_k-u'\| \leq M_2p_kd +O(h^2),\; 
%\end{split}
\end{equation}
\end{Theorem}
%\begin{Proof}
\noindent {\bf Proof.}
The first above estimate is easily obtained if representing
\begin{equation*}
U_k(t_i)-u(t_i)= (u_k(t_i)-u(t_i))+(U_k(t_i)-u_k(t_i))
\end{equation*}
and using  the first estimate in Theorem \ref{thm1} and the first estimate in \eqref{eq:prop3b}. The remaining estimates are obtained in the same way. Thus, the theorem is proved.  $\blacksquare$
%\end{Proof}

\section{Examples}

\noindent {\bf Example 1.} Consider the problem \eqref{eq:O3} with
\begin{align*}
&k(x,t)= e^x \sin (\pi t), \quad (x,t)\in [0,1] \times [0,1],\\
&f(x,u(x),u'(x) , \int_0^1 k(x,t) u(t) dt) = u^2(x) \int _0^1 k(x,t)u(t) dt  + u(x) u'(x) \\
&-\frac{1}{2} e^x \sin ^2 (\pi x) + \pi ^4 \sin (\pi x) - \frac{\pi}{2} \sin (2 \pi x).
\end{align*}
In this case 
$$ f(x,u,v,z)=u^2z+uv-\frac{1}{2} e^x \sin ^2 (\pi x) + \pi ^4 \sin (\pi x) - \frac{\pi}{2} \sin (2 \pi x)
$$
and $M_2  = \dfrac{2 e}{\pi} $.
It is possible to verify that the function $u= \sin (\pi x)$ is exact solution of the problem.  In the domain $\mathcal{D}_M$  defined by 
\begin{align*}
\mathcal{D}_M = \{(x,u,v,z)\; | \; &0\leq x\leq 1,\; |u|\leq M_0 M, |v| \le M_1M,
 |z|\leq M_0M_2 M \}
\end{align*}
we have
\begin{align*}
|f(x,u,v,z)| \le M_0^3M_2 M^3 + M_0M_1 M^2 + (\pi ^4 + \frac{\pi}{2}+\frac{e}{2}).
\end{align*}
It is possible to verify that for $M=113$ all the conditions of Theorem \ref{thm1} are satisfied with
%\begin{align*}
$ L_0 = 12.2010, L_1= 1.4714, L_2=2.1649, q=0.2690$.
%\end{align*}
Therefore, the problem has a unique solution $u(x)$ satisfying the estimates $|u(x)| \le 1.4714 , |u'(x)| \le 4.7083
$. These theoretical estimates are somewhat greater than the exact estimates $|u(x)| \le 1 , |u'(x)| \le \pi $.\par
Below we report the numerical results by the discrete iterative method  \eqref{iter1d}-\eqref{iter3d}  for the problem. In Tables \ref{table:1} and \ref{table:2} $Error =\|U_m-u \|$, where $u$ is the exact solution of the problem.
\begin{table}[ht!]
\centering
\caption{The convergence in Example 1 for the stopping criterion $\|U_m-u \|\le h^2$}
\label{table:1}
\begin{tabular}{cccc}
\hline 
$N$ &	$h^2$ &	$m$	&$Error$ \\
\hline 
50	& 4.0000e-04 &	2 &	1.4305e-04\\
100	& 1.0000e-04 &	3&	2.8588e-06\\
150	&4.4444e-05	&3	&2.8599e-06\\
200	&2.5000e-05	&3	&2.8602e-06\\
300	&1.1111e-05	&3	&2.8603e-06\\
400&	6.2500e-06&	3	&2.8603e-06\\
500	&4.0000e-06	&3	&2.8603e-06\\
800	&1.5625e-06&	4&	5.7485e-08\\
1000&	1.0000e-06	&4	&5.7486e-08\\
\hline 
\end{tabular} 
\end{table}
 It is interesting to notice that if taking the stopping criterion $\|\Phi_m-\Phi_{m-1} \|\le 10^{-10}$ then we obtain better accuracy of the approximate solution with more iterations. See Table \ref{table:2}. 

\begin{table}[ht!]
\centering
\caption{The convergence in Example 1 for the stopping criterion $\|\Phi_m-\Phi_{m-1} \|\le 10^{-10}$}
\label{table:2}
\begin{tabular}{cccc}
\hline 
$N$ &	$h^2$ &	$m$	&$Error$ \\
\hline 

50	& 4.0000e-04 &	7 &	2.2152e-08\\
100	& 1.0000e-04 &	7&	1.3831e-09\\
150	&4.4444e-05	&7	&2.7279e-10\\
200	&2.5000e-05	&7	&8.5995e-11\\
300	&1.1111e-05	&7	&1.6618e-11\\
400&	6.2500e-06&	7	&4.9447e-12\\
500	&4.0000e-06	&7	&1.7567e-12\\
800	&1.5625e-06&	7&	1.4588e-13\\
1000&	1.0000e-06	&7	&3.3318e-13\\
\hline 
\end{tabular} 
\end{table}
From Table \ref{table:2} we see that the accuracy of the approximate solution is near $O(h^4)$ although by the proved theory it is only $O(h^2)$.\\

\noindent {\bf Example 2.} (Example 4.2 in \cite{Wang}) Consider the nonlinear fourth order BVP
\begin{equation}\label{eq:E1}
\begin{split}
u^{(4)}(x)& =\sin (\pi x)[(2-u^2(x)) \int _0^1 tu(t) dt +1], x\in (0,1)\\
u(0)&=0 , \; u(1)=0, \; u''(0)=0 , \; u''(1)=0,
\end{split}
\end{equation}
This is the problem \eqref{eq:O3} with
\begin{align*}
&k(x,t)=\sin (\pi x) t, \quad (x,t)\in [0,1] \times [0,1],\\
&f(x,u(x), u'(x), \int_0^1 k(x,t) u(t) dt) = (2-u^2(x)) \int _0^1 \sin (\pi x)tu(t) dt  +\sin (\pi x).
\end{align*} 
So, $f(x,u,v,z) = (2-u^2)z+ \sin (\pi x)$. \\
It is easy to see that $M_2 =\max_{0\leq x\leq 1} \int_0^1 |k(x,t)|dt = \frac{1}{2}$.
Since $M_0$ and $ M_1$ are given by \eqref{eq:M} we define
\begin{align}\label{eq:A11}
\mathcal{D}_M = \{(x,u,v,z)\; | \; &0\leq x\leq 1,\; |u|\leq \frac{5}{384} M, |v| \le \frac{1}{24}M,
 |z|\leq \frac{5}{768} M \}.
\end{align}
It is possible to verify that for $M=1.1$ all the assumptions of Theorem \ref{thm1} are satisfied with
$ L_0=2.0515e-04, L_1=0, L_2=2, q= 0.0130$.
Therefore, the problem \eqref{eq:E1} has a unique solution satisfying $|u(x)| \le  0.0143,\; |u'(x)|\le 0.0458$.\par
{\it It is worth emphasizing that in \cite{Wang} by the monotone  method the author could only prove  the convergence of the iterative sequences to extremal solutions of the problem but not the existence and uniqueness of solution.}\par
Using the discrete iterative method \eqref{iter1d}-\eqref{iter3d} on the grid with grid step $h=0.01$ and the stopping criterion $\|\Phi_m-\Phi_{m-1} \|\le 10^{-10}$ we found an approximate solution after $7$ iterations. The graph of this approximate solution is depicted in Figure 1.
\begin{figure}[ht]
\begin{center}
\includegraphics[height=6cm,width=9cm]{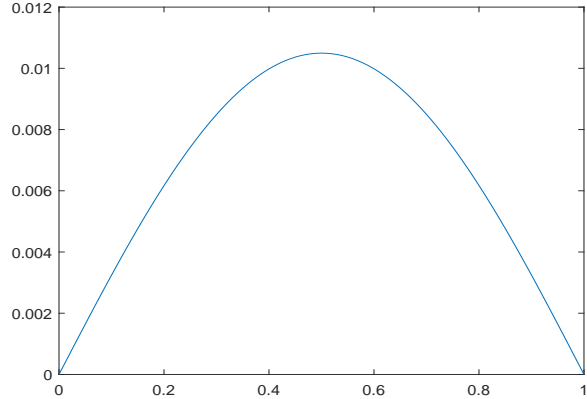}
\caption{The graph of the approximate solution in Example $2$. }
\label{fig1}
\end{center}
\end{figure}

\section{Conclusion}
In this paper we have established the existence and uniqueness of solution for a fourth order nonlinear integro-differential equation with the Navier boundary conditions and proposed an iterative method at both continuous and discrete levels for finding the solution. The second order of accuracy of the discrete method has been proved. Some examples, where the exact solution is known and is not known, demonstrate the validity of the obtained theoretical results and the efficiency of the iterative method. It should be emphasized that for the example of Wang in \cite{Wang} we have established the existence and uniqueness of solution and found it numerically but Wang could prove only the convergence of the iterative sequences constructed by the monotone method to extremal solutions.\par 
The method used in this paper with appropriate modifications can be applied to nonlinear integro-differential equations of any order with other boundary conditions and more complicated nonlinear terms. This is the direction of our research in the future.

\end{document}